\documentclass[11pt]{elsarticle}
\usepackage[T1]{fontenc}
\usepackage[utf8]{inputenc}
\usepackage{amsmath,amssymb,amsthm,mathtools,booktabs,array}
\usepackage[margin=1in]{geometry}
\usepackage[hidelinks]{hyperref}
\usepackage{enumitem}
\numberwithin{equation}{section}
\newtheorem{theorem}{Theorem}[section]
\newtheorem{proposition}[theorem]{Proposition}
\newtheorem{lemma}[theorem]{Lemma}

\theoremstyle{definition}
\newtheorem{definition}[theorem]{Definition}
\theoremstyle{remark}
\newtheorem{remark}[theorem]{Remark}
\newtheorem{example}[theorem]{Example}
\newcommand{\R}{\mathbb R}
\newcommand{\N}{\mathbb N}
\newcommand{\Nzero}{\mathbb N_0}
\newcommand{\E}{\mathbb E}
\newcommand{\Var}{\operatorname{Var}}
\newcommand{\rad}{\mathcal R}
\newcommand{\fall}[2]{#1^{\underline{#2}}}
\newcommand{\rise}[2]{(#1)_{#2}}

\begin{document}
\begin{frontmatter}
\title{Maximal Absolute Monotonicity of Balanced Rational Functions with a Single Repeated Pole}
\author{Rachid Ait-Haddou\corref{cor1}}
\ead{raithaddou@pmu.edu.sa}
\address{Department of Mathematics and Natural Sciences, College of Sciences and Human Studies, Prince Mohammad Bin Fahd University, Al Khobar 31952, Saudi Arabia}
\begin{abstract}
We study the radius of absolute monotonicity of rational approximations
to the exponential of the form
\[
\phi(z)=\frac{P(z)}{(1-z/b)^s},
\qquad b>0,\qquad \deg P\le s,
\]
with approximation orders two and three. 
For second-order approximation, we prove that the maximal radius in
degree $s=5$ is $10$ and determine the unique extremizer. Hence the
conjectured value $2s$ by L\'oczi and Ketcheson (LMS Journal of Computation and Mathematics, 17 (2014), no.~1, 159--205) remains valid in this case.
In contrast, for
every integer $s\ge6$ we construct explicitly a rational function with
rational coefficients and exact radius
\[
2s+\frac{4(s-3)}{s^2-3s+4}>2s.
\]
For third-order approximation, we construct, for every $s\ge5$, an
explicit rational family with exact radius
\[
2s-1>s-1+\sqrt{s^2-1}.
\]
These constructions disprove the conjectures 
of L\'oczi and Ketcheson concerning rational functions with a single
repeated pole. 	
\end{abstract}

\end{frontmatter}

\section{Introduction and main results}\label{sec:intro}
A rational approximation to the exponential may inherit important
monotonicity properties from the exact semigroup.  More precisely, if
a rational function $\phi$ approximates $e^z$ near the origin and is
absolutely monotonic on an interval $[-R,0]$, then $R$ provides a
step-size threshold under which positivity, monotonicity, or
contractivity properties of the forward Euler discretization can be
preserved in the corresponding linear problem.  This point of view goes
back to the classical theory of absolute monotonicity for numerical
time discretizations and, in the rational setting, to the work of
van de Griend and Kraaijevanger~\cite{VGK}.  In this sense, the radius
of absolute monotonicity is a natural scalar measure of the
monotonicity-preserving potential of a rational approximation to the
exponential.

For Runge--Kutta methods, however, one has to distinguish carefully
between the scalar stability function and the full method.  An
$s$-stage Runge--Kutta method with coefficient matrix $A$ and weight
vector $b$ has stability function
\[
\phi(z)=1+zb^T(I-zA)^{-1}\mathbf 1.
\]
The radius of absolute monotonicity of $\phi$ controls the corresponding
linear monotonicity problem, whereas the nonlinear SSP coefficient of
the method also depends on the internal stages.  In particular, the
latter requires positivity conditions on the full Runge--Kutta
realization and can be strictly smaller than the radius of its scalar
stability function.  L\'oczi and Ketcheson record this distinction in
\cite{LK}, where the method coefficient is bounded
above by the radius of absolute monotonicity of the associated
stability function.  The distinction is also emphasized in the
Runge--Kutta literature: the scalar stability function carries only the
linear information and does not encode the positivity constraints on
the internal stages.  In the present paper we study only this scalar
rational approximation problem.  In particular, we impose neither
A-stability nor the existence of a Runge--Kutta realization with
nonnegative internal coefficients.

The systematic study of absolute monotonicity for rational
approximations to the exponential was initiated by van de Griend and
Kraaijevanger~\cite{VGK}.  They investigated how large the radius of
absolute monotonicity can be when the numerator and denominator degrees
and the order of approximation are prescribed.  This problem is closely
related to monotonicity-preserving implicit one-step methods and
provides natural upper bounds for SSP Runge--Kutta methods.  L\'oczi
and Ketcheson~\cite{LK} returned to this problem and obtained a number
of sharp results for structured classes of rational functions.  Among
the classes considered there are rational functions for which the
denominator has a single real pole of prescribed multiplicity.  Such
functions arise naturally from singly diagonally implicit
Runge--Kutta methods, but the scalar optimization problem is meaningful
independently of whether a Runge--Kutta realization exists.

The balanced case, in which both numerator and denominator have degree
at most $s$, is especially simple to state.  In the single-pole class
one may write
\[
\phi(z)=\frac{P(z)}{(1-z/b)^s},
\qquad
b>0,\qquad \deg P\le s.
\]
For approximation order two, L\'oczi and Ketcheson conjectured that the
maximal radius in this class is $2s$; for approximation order three,
they proposed the value
\[
s-1+\sqrt{s^2-1}.
\]
They prove the validity of their proposed values for low degree cases $s \leq 4$. Their work also shows how questions of this
type are connected with the optimization of implicit SSP methods, but
the conjectures themselves are scalar statements about rational
functions. It is precisely these repeated-pole scalar conjectures that
we revisit in this work.

Our starting point comes from a different but closely related
quadrature interpretation. Ait-Haddou~\cite{AH} showed that the
optimization of polynomial maximal radius of monotonicity  can be reformulated as a positive quadrature problem with integer nodes for Poisson--Charlier measures. Thus a polynomial stability function can be viewed as the probability generating function of a finitely supported distribution on the nonnegative integers whose factorial moments agree with those of a Poisson law up to the required order.  This viewpoint is especially effective because positivity of the coefficients is built into the
measure-theoretic formulation.  The Poisson--Charlier interpretation is
one of the central ingredients in the polynomial theory developed
in~\cite{AH}; in particular, that work shows that positive quadrature
with integer nodes is equivalent to the absolute-monotonicity
representation underlying optimal linear SSP polynomials.

The rational repeated-pole problem suggests a natural extension of this
idea.  Instead of restricting the positive measure to have finite
support, we allow an infinite distribution on the nonnegative integers,
but require its probability generating function to be rational. If
\[
F(w)=\phi(r(w-1))=\sum_{k=0}^{\infty}\pi_k w^k,
\qquad \pi_k\ge0,\qquad \sum_{k\ge0}\pi_k=1,
\]
then approximation of $e^z$ to order $p$ is equivalent to the factorial
moment conditions
\[
\sum_{k=0}^{\infty} \frac{k!}{(k-j)!}\pi_k=r^j,
\qquad 0\le j\le p.
\]
Thus the scalar rational problem becomes a positive moment problem on
the integer lattice.  In the repeated-pole case, the tail of the
sequence $(\pi_k)$ is a polynomial modification of a geometric, or more
generally negative-binomial, distribution.  This is the point at which
Meixner-type structure enters naturally.  It provides both a source of
explicit constructions and a convenient mechanism for proving
coefficientwise positivity.

This extension leads to a phenomenon not visible in the previously
verified low-degree cases. We show that the conjectured second-order value $2s$
remains correct through degree five, but starting at degree six one can
construct positive rational generating functions whose radii exceed
$2s$.  At third order, the proposed repeated-pole value already fails
in degree five. The counterexamples are not obtained by numerical
optimization: they arise from explicit positive modifications of
negative-binomial laws, and their radii are proved exactly by analyzing
the first nonzero coefficient of the shifted generating function.
Accordingly, the present work should be viewed as an extension of the
positive-quadrature viewpoint of Ait-Haddou \cite{AH} from finitely supported polynomial
generating functions to rational generating functions with structured
infinite tails.

\begin{definition}
For a real rational function $\phi$ analytic on the nonpositive real axis, we define the {\it{radius of absolute monotonicity}} $\rad(\phi)$ of $\phi$ as
\[
\rad(\phi)=\sup\bigl(\{0\}\cup\{r>0:\phi^{(j)}(x)\ge0\text{ for all }j\in\Nzero\text{ and }x\in[-r,0]\}\bigr).
\]
For $s\in\N$ and $p\in \Nzero$, let
\begin{equation*}
\mathcal{B}_{s,p}=
\left\{
\frac{P(z)}{(1-z/b)^s}, b>0,\ P\in\R[z],\ \deg P\le s,\ \phi(z)-e^z=O(z^{p+1}) \quad\text{as }z\to0
\right\}.
\end{equation*}
The {\it{optimal radius of absolute monotonicity}} $B_{s,p}$ is defined by 
\[
B_{s,p}=\sup\bigl(\{0\}\cup\{\rad(\phi):\phi\in\mathcal B_{s,p}\}\bigr).
\]
\end{definition}

Our main results in this work as summarised as follows. 

\begin{theorem}\label{thm:main}
The following statements hold.
\begin{enumerate}[label=\textup{(\roman*)}]
\item $B_{5,2}=10$. The only member of $\mathcal B_{5,2}$ absolutely monotonic on $[-10,0]$ is
\begin{equation}\label{eq:mid5}
\phi_*(z)=\left(\frac{1+z/10}{1-z/10}\right)^5.
\end{equation}
Its radius of absolute monotonicity is 10.
\item For every integer $s\ge6$,
\begin{equation}\label{eq:order2bound}
B_{s,2}\ge 2s+\frac{4(s-3)}{s^2-3s+4}>2s.
\end{equation}
For each such $s$, an explicitly specified function with rational coefficients has exactly the radius on the right.
\item For every $s\in\N$,
\begin{equation}\label{eq:productbound}
B_{s,2}\ge 2s+\left\lfloor\frac{s}{7}\right\rfloor.
\end{equation}
\item For every integer $s\ge5$,
\begin{equation}\label{eq:order3bound}
B_{s,3}\ge 2s-1>s-1+\sqrt{s^2-1}.
\end{equation}
The lower value $2s-1$ is the exact radius of an explicit rational-coefficient function of exact order three. In particular, $9\le B_{5,3}\le10$.
\end{enumerate}
\end{theorem}
\begin{remark}
In \cite{LK}, the authors considered the family of rational
functions
\[
\hat{\Pi}_{s,p}
=
\left\{
\phi(z)=\frac{P(z)}{(1-az)^s}:
\begin{array}{l}
a\in\mathbb R,\quad P\in\mathbb R[z],\quad \deg P\le s,\\
\phi(z)-e^z=O(z^{p+1})\quad\text{as }z\to0
\end{array}
\right\}.
\]
If $a<0$, then $\phi$ has no positive real pole and,
by \cite{LK}, its radius of absolute monotonicity is zero.	
If $a=0$, then $\phi$ is a polynomial of degree at most $s$.
The classical upper bound for polynomials of order $p$ 
gives \cite{K}
\[
\rad(\phi)\le s-p+1.
\]
For $p=2,3$, these bounds are strictly smaller than the
corresponding radii established in
Theorem~\ref{thm:main}. Consequently, neither $a<0$
nor $a=0$ can yield an optimal function.
We may therefore restrict our attention, without loss
of generality, to $a>0$, or equivalently to the family
$\mathcal B_{s,p}$ obtained by setting $b=1/a>0$.
\end{remark}

\medskip
\noindent\textbf{Notation.}
Throughout the paper, we use the falling and rising factorials
\[
\fall{x}{k}=\prod_{j=0}^{k-1}(x-j),
\qquad
\rise{x}{k}=\prod_{j=0}^{k-1}(x+j),
\]
with the convention that empty products equal $1$.
Removable singularities of rational functions are understood
to be cancelled when considering analyticity or derivatives.

For a power series
\[
F(w)=\sum_{i=0}^{\infty}p_iw^i,
\qquad p_i\ge0,\qquad F(1)=1,
\]
convergent in a neighborhood of $w=1$, we denote by $X_F$
the nonnegative integer-valued random variable defined by
\[
\mathbb P(X_F=i)=p_i=[w^i]F.
\]
Thus $F(w)=\E[w^{X_F}]$, and its factorial moments satisfy
\[
\E[\fall{X_F}{k}]=F^{(k)}(1).
\]
In particular,
\[
\E[X_F]=F'(1),
\qquad
\Var[X_F]=F''(1)+F'(1)-F'(1)^2.
\]

\section{Positive generating functions and exact radii of absolute monotonicity}
\label{sec:pgf}
We start with some preliminary lemmas.
\begin{lemma}\label{lem:pgf}
Let $b,r>0$, $s\in\N$, and 
\[
\phi(z)=\frac{P(z)}{(1-z/b)^s} \quad \text {with} \quad  \phi(0)=1.
\]
Set
\[
\eta=\frac r{b+r},\qquad F(w)=\phi(r(w-1)).
\]
Then $\phi$ is absolutely monotonic on $[-r,0]$ if and only if all Taylor coefficients
$\pi_k=[w^k]F$ are nonnegative. In that case they sum to one, and for the random variable $X_F$ with $\mathbb P(X_F=k)=\pi_k$,
\begin{equation}\label{eq:factorial}
\phi^{(j)}(0)=1\ (0\le j\le p)
\quad\Longleftrightarrow\quad
\E\fall {X_F} j=r^j\ (0\le j\le p).
\end{equation}
In particular, order at least two is equivalent to $\E[X_F]=\Var[X_F]=r$.
\end{lemma}
\begin{proof}
The only possible pole of $F$ is $\eta^{-1}>1$. Its Taylor series therefore converges absolutely on a disk containing $[0,1]$, and can be differentiated termwise there to any finite order. Its coefficients are $\pi_k=r^k\phi^{(k)}(-r)/k!$. Nonnegativity of the coefficients is necessary. It is sufficient because each differentiated series is a sum of nonnegative terms on $[0,1]$. Evaluation at one gives 
\[
\sum_{k=0}^{\infty}\pi_k=F(1)=\phi(0) = 1,
\] 
and differentiation at one gives $\E[\fall {X_F} j]=F^{(j)}(1)=r^j\phi^{(j)}(0)$. All moments exist because the series converges on a disk of radius greater than one. Finally, $\E[X_F^2]=\E[\fall {X_F}2]+\E[X_F]=r^2+r$ is equivalent to the claimed variance identity.
\end{proof}

\begin{lemma}
\label{lem:gap}
Let $F$ have nonnegative Taylor coefficients, radius of convergence greater than one, and $F(1)=1$. If $d\in\Nzero$ satisfies $[w^d]F=0$ and $[w^{d+1}]F>0$, then $\phi(z)=F(1+z/r)$ has $\rad(\phi)=r$ for every $r>0$ for which it is considered. The same conclusion holds if $F(w)=c w^q+O(w^{q+1})$ with $c>0$ and $q\ge1$.
\end{lemma}
\begin{proof}
The nonnegative series argument of Lemma~\ref{lem:pgf} proves the bound $R(\phi) \ge r$. In the first case,
$F^{(d)}(w)=(d+1)![w^{d+1}]F\,w+O(w^2)$, which is negative immediately to the left of zero. In the second case use $d=q-1$. The change of variables $F(w)=\phi(r(w-1))$  carries these points to points immediately to the left of $-r$ and multiplies derivatives by positive powers of $r$.This concludes the proof.
\end{proof}
\begin{lemma}\label{lem:mixture}
	Let $F(w)=\sum_{i=0}^{\infty} a_i w^i, \; a_i\geq 0,$
	be a power series defined in a neighborhood of \(w=1\) with $F(1) = 1$. Set
	\[
	M:=F'(1), \quad  
	V:=F''(1)+F'(1)-F'(1)^2.
	\]
	Assume that $M>1$ and $0<V\leq M$.
	Define
	\[
	\epsilon:=\frac{M-V}{M^2},
	\quad
	G(w):=\epsilon+(1-\epsilon)F(w),
	\quad 
	r:=(1-\epsilon)M
	=M-1+\frac{V}{M}.
	\]
	Then \(G(1)=1\), all the coefficients of \(G\) are nonnegative, and
	\[
	G'(1)=r, \quad G''(1)+G'(1)-G'(1)^2=r.
	\]
\end{lemma}
\begin{proof} Since $M>1$ and $0<V\leq M$,
	we have
	\[
	0\leq \epsilon
	=\frac{M-V}{M^2}
	<\frac{1}{M}<1.
	\]
	Hence $G(w)=\epsilon+(1-\epsilon)F(w)$ has nonnegative coefficients, and
	$G(1)=\epsilon+(1-\epsilon)F(1)=1$. Moreover,
	\[
	G'(1)=(1-\epsilon)F'(1)
	=(1-\epsilon)M=r.
	\]
	Since $ F''(1)=V-M+M^2$, we also have
	\[
	G''(1)
	=(1-\epsilon)F''(1)
	=(1-\epsilon)(V-M+M^2).
	\]
	Therefore
	\[
	G''(1)+G'(1)-G'(1)^2
	=(1-\epsilon)(V-M+M^2) +(1-\epsilon)M
	-(1-\epsilon)^2M^2
	=(1-\epsilon)
	\bigl(V+\epsilon M^2\bigr).
	\]
	By the definition of \(\epsilon\), $\epsilon M^2=M-V$,
	and hence
	\[
	G''(1)+G'(1)-G'(1)^2=(1-\epsilon)M=r.
	\]
	This proves the result.
\end{proof}
\begin{remark}
The previous lemma can be re-casted completely in terms of the probability associated with $F$ and $G$ as follows: 
Let $X_{F}$ be a nonnegative integer-valued random variable with $\E[X_F]>1$ and $0<\Var[X_F]\le \E[X_F] $. Set
\[
\epsilon=\frac{\E[X]-\Var[X]}{\E[X]^2},\qquad r=(1-\epsilon)=\E[X]-1+ \frac{\Var[X]}{\E[X]}.
\]
Let \(X_G\) be a random variable with probability law 
$\mathcal L(X_G) = \epsilon\delta_0+(1-\epsilon)\mathcal L(X_F)$. 
Then $X_G$ has mean and variance $r$. Here $\mathcal L(X)$ denotes the law of $X$.
\end{remark}

\section{Second-order counterexamples}\label{sec:second}
In our construction of counter-examples, we shall use the negative-binomial distribution with shape $q>0$ and parameter $0<\eta<1$, i.e.;
\begin{equation}\label{eq:nb}
\mathbb P(V=k)=(1-\eta)^q\frac{\rise qk}{k!}\eta^k,
\qquad
\E[w^V]=\left(\frac{1-\eta}{1-\eta w}\right)^q.
\end{equation}
The binomial series proves both the normalization and the generating formula in (\ref{eq:nb}). 
Differentiation at one gives
\begin{equation}\label{eq:ffactorial}
\E[\fall Vj]=\rise qj\left(\frac{\eta}{1-\eta}\right)^j.
\end{equation}
\subsection{A sevenfold-pole counter-example}
Let $V$ have the negative-binomial distribution with shape $q=5$ and parameter $\eta = 1/2$. Thus, from (\ref{eq:nb})
\[
\mathbb P (V=k)=\binom{k+4}{4}2^{-k-5},
\qquad k=0,1,2,\ldots,
\]
and its probability generating function is
\[
G(w)=\mathbb E[w^V]=(2-w)^{-5},
\qquad |w|<2.
\]
Since $(V-2)(V-3)=V(V-1)-4V+6$, differentiation gives
\[
H(w)
:=\mathbb E[(V-2)(V-3)w^V]=w^2G''(w)-4wG'(w)+6G(w)=\frac{8(7w^2-8w+3)}{(2-w)^7}.
\]
In particular, $H(1)=16$. Therefore
\begin{equation}\label{F7generating}
F_7(w):=\frac{w^5H(w)}{16}
=\frac{w^5(7w^2-8w+3)}{2(2-w)^7}
\end{equation}
has the expansion
\[
F_7(w)=
\sum_{k=0}^{\infty}
\frac{(k-2)(k-3)}{16}
\binom{k+4}{4}2^{-k-5}w^{k+5}.
\]
Every coefficient is nonnegative because
$(k-2)(k-3)\ge0$ for every nonnegative integer $k$.
Moreover, $F_7(1)=H(1)/16=1$.
Thus $F_7$ is obtained by reweighting the probabilities of $V$,
normalizing them, and shifting their locations by five.
Define
\begin{equation}
\label{subphi}
\phi_7(z)=F_7(1+z/15).
\end{equation}
Exact differentiation give  
\[
\phi_7(0)=\phi_7'(0)=\phi_7''(0)=1,
\qquad
\phi_7'''(0)=\frac{226}{225}\ne1.
\]
 Hence $\phi_7 \in \mathcal{B}_{7,2}$. Since $F_7(w)=3w^5/256+O(w^6)$, 
 Lemma~\ref{lem:gap} gives exact radius $\rad(\phi_7)=15 > 14$. This already provide a counter-example to
 L\'oczi and Ketcheson conjecture. Substitution of $F_7$ in (\ref{subphi}) yields
\begin{equation}\label{eq:sevenphi}
\phi_7(z)=
\frac{(1+z/15)^5(1+z/5+7z^2/450)}{(1-z/15)^7}.
\end{equation}


\subsection{A rational family in every degree at least six}
We shall now generalise the previous construction to any $s \geq 6$. 
For integer $s\ge6$, define
\begin{gather}
D_s=s^2-5s+8,\qquad
r_s=2s+\frac{4(s-3)}{s^2-3s+4},\label{eq:secondr}\\
\epsilon_s=\frac{2(s^2-7s+8)}{(s-2)(s^2-3s+4)^2},\label{eq:secondeps}\\
H_s(w)=\frac{w^{s-2}[s(s-1)w^2-4sw+8]}{D_s(2-w)^s},\quad
F_s(w)=\epsilon_s+(1-\epsilon_s)H_s(w).
\label{eq:secondfamily}
\end{gather}
Define 
\begin{equation}\label{order2degrees}
\phi_s(z)=F_s(1+z/r_s).
\end{equation}
\begin{proposition}\label{prop:secondfamily}
The function in \eqref{order2degrees} belongs to $\mathcal B_{s,2}$ and has exact radius of absolute monotonicity $r_s$. 
Its numerator and denominator have degree exactly $s$, with no pole cancellation.
\end{proposition}
\begin{proof}
Let $V$ have the negative-binomial distribution with shape $q = s-2$ and $\eta=1/2$ in \eqref{eq:nb}.
Its generating function is given by 
\[
G_s(w)=\mathbb E[w^V]=(2-w)^{-s+2},
\qquad |w|<2.
\] 
Since $(V-1)(V-2) = V(V-1)-2 V + 2$, the factorial-moment formula (\ref{eq:ffactorial}) and differentiation 
show that the quantities $D_s$ and $H_s$ defined in (\ref{eq:secondr}) and (\ref{eq:secondfamily}) are given by 
\[
D_s = \E[(V-1)(V-2)] \quad \text{and} 
\quad
H_s(w)=\frac{w^{s-2}}{D_s}\E[(V-1)(V-2)w^V].
\]
Thus $H_s$ is a probability generating function with distribution $X_s$ given by
\begin{equation*}
\mathbb P(X_s=s-2+k)=\frac{(k-1)(k-2)}{D_s}\; \mathbb P(V=k) = \frac{(k-1)(k-2)}{D_s} \binom{k+s-3}{s-1}2^{-k-s+2}, 
\; k=0,1,2,\ldots,
\end{equation*} 
The mean of $X_s$ is given by 
\begin{equation*}
\begin{split}
\mathbb{E}[X_s] & = \frac{1}{D_s}\sum_{k=0}^{\infty} (k+s-2) (k-1)(k-2) \mathbb{P}(V = k) \\
& = \frac{1}{D_s} \left( (s-2) \sum_{k=0}^{\infty} (k-1)(k-2) \mathbb{P}(V = k) + \sum_{k=0}^{\infty} \fall k3 \mathbb{P}(V = k) \right) \\  
& = \frac{1}{D_s} (s-2) \E[(V-1)(V-2)] + \frac{1}{D_s} \E[\fall V3] 
= (s-2) + \frac{\rise{s-2}3}{D_s}.  
\end{split}
\end{equation*}
Thus, we have 
\begin{equation}
\mathbb{E}[X_s] =  (s-2) + \frac{\rise{s-2}3}{D_s}  = \frac{2(s-2)(s^2-3s+4)}{D_s}.\label{eq:secondM}
\end{equation}
Similarly, using that $V^2(V-1)(V-2) = \fall V4 + 3 \fall V3$, we have 
\begin{equation*}
\begin{split}
\mathbb{E}[X^2_s] & = \frac{1}{D_s}\sum_{k=0}^{\infty} (k+s-2)^2 (k-1)(k-2) \mathbb{P}(V = k) \\
	& = \frac{1}{D_s} \left( 2 (s-2) \E[\fall V3]+ \E[\fall V4 + 3 \fall V3] + (s-2)^2 \E[(V-1)(V-2)] \right) \\ 
	& = \frac{2(s-2)(s-2)_3 + (s-2)_4+ 3(s-2)_3 + (s-2)^2D_s}{D_s} \\ 
	& = (s-2)^2 +  \frac{3 s \rise{s-2}3}{D_s}.  
\end{split}  
\end{equation*}
Thus, the variance of $X_s$ is   
\begin{equation}
V[X_s]:= \mathbb{E}[X^2_s] - \mathbb{E}[X_s]^2 = \frac{2s(s-2)(s-1)(s^2-7s+16)}{D_s^2}. \label{eq:secondV}
\end{equation} 
Subtraction gives
\[
\mathbb E[X_s]-V[X_s]=\frac{8(s-2)(s^2-7s+8)}{D_s^2}>0 \quad \text{for} \quad  s \geq 6.
\]
Clearly from (\ref{eq:secondM}), we have $\mathbb E[X_s] >1$, since $\mathbb E[X_s] \geq  s-2 \geq 4$. Moreover, the distribution of $X_s$ is nonconstant, so $V[X_s]>0$. The value $\epsilon_s$ given in (\ref{eq:secondeps}) equals $(E[X_s]-V[X_s])/E[X_s]^2$. Lemma~\ref{lem:mixture} therefore gives mean and variance $r_s = (1- \epsilon_s)\mathbb E[X_s]$ for $F_s$. This quantity $r_s$ coincides with the one given in (\ref{eq:secondr}). Lemma~\ref{lem:pgf} proves second-order accuracy and absolute monotonicity over $[-r_s,0]$.
Its first nonconstant coefficient is that of $w^{s-2}$ and is positive. Since $s-2\ge4$, it is preceded by a zero coefficient. Lemma~\ref{lem:gap} proves exact radius of absolute monotonicity $\rad(\phi_s) = r_s$. At $w=2$ the quadratic in $H_s$ equals $4(s-1)(s-2)>0$, so the pole is genuine. The leading numerator coefficient over the denominator $D_s(2-w)^s$ is
\[
\epsilon_sD_s(-1)^s+(1-\epsilon_s)s(s-1).
\]
For even $s$ it is positive. For odd $s$, use $D_s<s(s-1)$ and $\epsilon_s<1/E[X_s]<1/2$ to obtain positivity. Thus the numerator degree is exactly $s$. This concludes the proof.
\end{proof}
\begin{remark}
The inequality $r_s > 2s$ in the previous proposition follows directly from \eqref{eq:secondr}. 
This proves part~(ii) of Theorem~\ref{thm:main}.
\end{remark}
\begin{example}
Take $s = 6$, we obtain 
\[
r_6=\frac{138}{11} > 12,\qquad
F_6(w)=\frac1{484}+\frac{69}{484}\frac{w^4(15w^2-12w+4)}{(2-w)^6}.
\]
The corresponding rational function in $\mathcal B_{s,2}$ is
\[
\phi_6(z)=F_6(1+z/r_6).
\]
\end{example}

\subsection{Product construction and asymptotic bound}
Let $F_1(w)=w/(2-w)$. $F_1$ is the generating function of the discrete probability $X_{F_1}$ 
whose mean and variance are equal two. Write $s=7q+d$ with $q=\lfloor s/7\rfloor$ and $0\le d\le6$. Let $F_7$ be the generating function defined in (\ref{F7generating}). The product
\[
F(w)=F_7(w)^qF_1(w)^d
\]
is a probability generating function of a sum of independent variables. 
Its mean and variance are $15q+2d=2s+q$. It has numerator and denominator degree $s$, no pole cancellation at $w=2$, and a positive zero at zero. Lemmas~\ref{lem:pgf} and~\ref{lem:gap} prove that the function 
\begin{equation*}
\phi(z) = F\left(1 + \frac{z}{2s + q}\right)
\end{equation*}
belong to $\mathcal B_{s,2}$ and $\rad(\phi) = 2s +q$. This proves part part~(iii) of Theorem~\ref{thm:main}. 
In particular, we have
\[
\liminf_{s\to\infty}\frac{B_{s,2}}s\ge\frac{15}{7}.
\]
This is a constructive lower bound, not a growth law for the optimum.
\section{Third-order counterexamples}\label{sec:third}
For integer $s\ge5$, set
\begin{gather}
u_s=\frac{s-1}{s},\qquad r_s=2s-1,\qquad b_s=\frac{s(2s-1)}{s-1},\label{eq:thirdparams}\\
c_s=\frac3{2(s-2)},\quad
d_s=-\frac{8s^2-17s+5}{(s-2)(s-1)(2s-1)},\quad
a_s=\frac{8s^3-4s^2-26s+7}{2(s-2)(2s-1)^2},\label{eq:thirdco}\\
\epsilon_s=\frac{2(s^2-4s+1)}{(s-2)(s-1)(2s-1)^2},\qquad
F_s(w)=\epsilon_s+\frac{w^{s-2}(c_s+d_sw+a_sw^2)}{(1+u_s-u_sw)^s}.\label{eq:thirdF}
\end{gather}
Define 
\begin{equation}
\label{order3counter}
\phi_s(z)=F_s(1+z/r_s).
\end{equation}

\begin{proposition}\label{prop:thirdfamily}
For every $s\ge5$, the function $\phi_s$ in (\ref{order3counter}) belongs $\mathcal B_{s,3}$  with $\rad(\phi_s) = 2s -1$.
Its numerator and denominator have degree exactly $s$, with no pole cancellation. 
\end{proposition}
\begin{proof}
For nonnegative integers $k,d$, define
\[
D_{s,d,k}(u)=\sum_{\ell=0}^{k}\binom{k}{\ell}\; \fall d\ell \; \rise{s}{k-\ell} \; u^{k-\ell}.
\]
The product rule identifies this expression as the $k$th derivative at $w=1$ of $w^d(1+u-uw)^{-s}$.
Consequently, $F_s(1) = \epsilon_s+c_s+d_s+a_s$ and for $k\ge1$,
\[
F_s^{(k)}(1) = c_s D_{s,s-2,k}(u_s) + d_s D_{s,s-1,k}(u_s) +a_s D_{s,s,k}(u_s).
\]
Substituting \eqref{eq:thirdparams}--\eqref{eq:thirdF}, clearing their positive denominators, and collecting powers of $s$ gives
\begin{align}
\epsilon_s+c_s+d_s+a_s&=1,\label{eq:thirdders}\\
c_sD_{s,s-2,k}(u_s)+d_sD_{s,s-1,k}(u_s)+a_sD_{s,s,k}(u_s)&=(2s-1)^k,
\qquad k=1,2,3.\nonumber
\end{align}
Moreover, we have
\begin{equation}\label{eq:thirddefect}
F_s^{(4)}(1)-r_s^4=
\frac{(2s-1)(4s^3-6s^2-12s+1)}{s^2}>0,
\end{equation}
since for $s\ge5$ the cubic factor is positive, i.e;  
\[
4s^3-6s^2-12s+1=4(s-5)^3+54(s-5)^2+228(s-5)+291 > 0.
\]
Equations~\eqref{eq:thirdders}--\eqref{eq:thirddefect} prove exact order three. 
Now we shall prove positivity of the Taylor coefficients of $F_s$. Set 
\[
\eta_s=\frac{u_s}{1+u_s}
=\frac{s-1}{2s-1},
\qquad
1-\eta_s= \frac{1}{1 + u_s} = \frac{s}{2s-1}.
\]
Then
\[
F_s(w)=
\epsilon_s+
(1-\eta_s)^s
\frac{w^{s-2}(c_s+d_sw+a_sw^2)}{(1-\eta_s w)^s}.
\]
The binomial expansion gives
\[
(1-\eta_s)^s (1- \eta_s w)^{-s}
=\sum_{m=0}^{\infty}B_mw^m,
\qquad
B_m=(1-\eta_s)^s
\binom{m+s-1}{s-1}\eta_s^m.
\]
Extracting the coefficient of
$w^{s-2+n}$ yields
\[
\pi_{s-2+n}:=[w^{s-2+n}]F_s
=c_sB_n+d_sB_{n-1}+a_sB_{n-2},
\qquad n\ge0,
\]
with the convention $B_m=0$ for $m<0$.
The constant term $\epsilon_s$ does not contribute because
$s-2+n\ge3$.
Introduce the positive weights
\[
\rho_s(n)=
\binom{n+s-3}{s-3}(1-\eta)^{s-2}\eta^n
=
\frac{(s-2)_n}{n!}
(1-\eta)^{s-2}\eta^n.
\]
Cancelling factorials in the binomial coefficients gives
\begin{align*}
	\frac{B_n}{\rho_s(n)}
	&=\frac{(1-\eta)^2}{(s-1)(s-2)}
	(n+s-2)(n+s-1),\\
	\frac{B_{n-1}}{\rho_s(n)}
	&=\frac{(1-\eta)^2}{(s-1)(s-2)}
	\frac{n(n+s-2)}{\eta},\\
	\frac{B_{n-2}}{\rho_s(n)}
	&=\frac{(1-\eta)^2}{(s-1)(s-2)}
	\frac{n(n-1)}{\eta^2}.
\end{align*}
Consequently,
\[
\pi_{s-2+n}
=
\rho_s(n)\frac{(1-\eta)^2}{(s-1)(s-2)}R_s(n),
\]
where
\[
R_s(n)=
c_s(n+s-2)(n+s-1)
+\frac{d_s}{\eta}n(n+s-2)
+\frac{a_s}{\eta^2}n(n-1).
\]
Using
\[
c_s=\frac3{2(s-2)},\qquad
\frac{d_s}{\eta}
=-\frac{8s^2-17s+5}{(s-2)(s-1)^2},
\qquad
\frac{a_s}{\eta^2}
=\frac{8s^3-4s^2-26s+7}{2(s-2)(s-1)^2},
\]
we find
\[
\begin{aligned}
	2(s-2)(s-1)^2R_s(n)
	=&3(s-1)^2(n+s-2)(n+s-1)\\
	&-2(8s^2-17s+5)n(n+s-2)\\
	&+(8s^3-4s^2-26s+7)n(n-1)\\
	=&(s-2)Q_s(n),
\end{aligned}
\]
where
\[
Q_s(n)=
s(8s-1)n^2-(18s^2-13s+2)n+3(s-1)^3.
\]
Thus $R_s(n)=Q_s(n)/(2(s-1)^2)$, and
\begin{equation}
\label{eq:thirdtails}
[w^{s-2+n}]F_s
=
\rho_s(n)
\frac{s^2Q_s(n)}
{2(s-2)(s-1)^3(2s-1)^2},
\qquad n\ge0.
\end{equation}
Every factor outside $Q_s(n)$ is strictly positive. To determine
its sign on the nonnegative integers, observe that
\[
Q_s(0)=3(s-1)^3>0,
\qquad
Q_s(1)=(s-1)(s-5)(3s-1)\ge0,
\]
and, for every integer $n\ge1$,
\[
Q_s(n+1)-Q_s(n)
=
6s^2+10s-2+2s(8s-1)(n-1)>0.
\]
Hence $Q_s(n)\ge0$ for every $n\ge0$, proving nonnegativity
of all the tail coefficients $[w^{s-2+n}]F_s, \; n \geq 0$.
Finally,
\[
[w^0]F_s=\epsilon_s
=\frac{2(s^2-4s+1)}
{(s-2)(s-1)(2s-1)^2}>0,
\qquad
[w^j]F_s=0\quad(1\le j\le s-3).
\]
Therefore every Taylor coefficient of $F_s$ is nonnegative and the coefficients at $1,\ldots,s-3$ vanish and that at $s-2$ is positive, 
proving that $\rad(\phi_s) = 2s-1$ by Lemma~\ref{lem:gap}.
The positive quadratic leading term in \eqref{eq:thirdtails} produces a genuine pole of order $s$. The leading numerator coefficient of \eqref{eq:thirdF} is $a_s+\epsilon_s(-u_s)^s$. Here $0<u_s<1$, $0<\epsilon_s<1$, and
\[
a_s-1=\frac{20s^2-44s+11}{2(s-2)(2s-1)^2}>0.
\]
It cannot vanish. This proves the degree assertion. This concludes the proof.
\end{proof}
\begin{remark}
Since  $\sqrt{s^2-1}<s$, Proposition~\ref{prop:thirdfamily} proves
part~(iv) of Theorem~\ref{thm:main}.
\end{remark}

\begin{example}
For $s = 5$, we have 
\begin{equation*}
F_5(w) = \frac{1}{81} - \frac{w^3 \left(\frac{259}{162} w^2 - \frac{10}{9} w + \frac{1}{2}\right) }{\left(\frac{4}{5} w - \frac{9}{5}\right)^5} 
\end{equation*}
The associated $\phi_5$ is given by
\begin{equation}\label{eq:thirdfive}
\phi_5(z)=\frac{1+\frac59z+\frac{109}{810}z^2+\frac{599}{36450}z^3
+\frac{1727}{1640250}z^4+\frac{9967}{369056250}z^5}
{(1-\frac4{45}z)^5}.
\end{equation}
According to the previous proposition, we have $\rad(\phi_5) = 9$. As an additional direct check, $F_5$ has coefficients
\[
\pi_0=\frac1{81},\quad \pi_1=\pi_2=0,\quad
\pi_{n+3}=\frac{3125}{118098}\left(\frac49\right)^n
\frac{(n-1)(n+1)(n+2)(65n-64)}{128},\quad n\ge0.
\]
Thus the taylor coefficient is positive at $n=0$, zero at $n=1$, and positive for $n\ge2$. Also
$\phi_5(z)-e^z=97z^4/145800+O(z^5)$. 
\end{example}

\section{The exact fifth-degree second-order optimum}\label{sec:s5}
This section proves part~(i) of the main Theorem~\ref{thm:main}.

\begin{lemma}\label{lem:s5tail}
Let $r>0$, and suppose that $\phi\in\mathcal B_{5,2}$ is
absolutely monotonic on $[-r,0]$. Its shifted generating function
\[
F(w)=\phi(r(w-1))
=\sum_{k=0}^{\infty}\pi_k w^k
\]
satisfies
\[
\pi_{k+1}=g_u(k)\Omega(k),
\qquad
g_u(k)=\frac1{1+u}\left(\frac{u}{1+u}\right)^k,
\qquad u=\frac rb>0,
\]
where $\Omega$ is a real polynomial of degree at most four and
$\Omega(k)\ge0$ for every $k\in\Nzero$.
\end{lemma}
\begin{proof}
Write
\[
\phi(z)=\frac{P(z)}{(1-z/b)^5},
\qquad b>0,\qquad \deg P\le5,
\]
and set
\[
u=\frac rb,\qquad
\eta=\frac{u}{1+u}\in(0,1).
\]
Since
\[
1-\frac{r(w-1)}b=(1+u)(1-\eta w),
\]
we have
\[
F(w)=\frac{A(w)}{(1-\eta w)^5},
\qquad
A(w)=\frac{P(r(w-1))}{(1+u)^5},
\qquad \deg A\le5.
\]
The Taylor series of $F$ converges for $|w|<\eta^{-1}$.
Moreover, absolute monotonicity gives $\pi_k \geq 0$ for $k \geq 0$. 
Because $\pi_0=F(0)=A(0)$, the polynomial $A(w)-\pi_0(1-\eta w)^5$
has zero constant term and degree at most five. It is therefore
divisible by $w$, so
\[
F(w)=\pi_0+\frac{wB(w)}{(1-\eta w)^5},
\qquad
B(w)=\sum_{j=0}^{4}\beta_jw^j.
\]
Thus
\[
\pi_{k+1}
=
\sum_{j=0}^{4}\beta_j
[w^k]\bigl(w^j(1-\eta w)^{-5}\bigr).
\]	
The binomial expansion
\[
(1-\eta w)^{-5}
=
\sum_{m=0}^{\infty}\binom{m+4}{4}\eta^m w^m
\]
shows that, for $k\ge j$,
\[
[w^k]\bigl(w^j(1-\eta w)^{-5}\bigr)
=
\eta^{k-j}\binom{k-j+4}{4}.
\]
Interpret the binomial coefficient as the polynomial
\[
\binom{x-j+4}{4}=
\frac{(x-j+1)(x-j+2)(x-j+3)(x-j+4)}{24}.
\]
If $0\le k<j\le4$, one of these four factors vanishes.
Hence the same coefficient formula remains valid for every
$k\ge0$, including those indices below $j$.	
Consequently,
\[
\pi_{k+1}
=
\eta^k\sum_{j=0}^{4}
\beta_j\eta^{-j}\binom{k-j+4}{4}.
\]
Define
\[
\Omega(x)=
\frac1{1-\eta}
\sum_{j=0}^{4}
\beta_j\eta^{-j}\binom{x-j+4}{4}.
\]
This is a real polynomial of degree at most four, and
\[
\pi_{k+1}
=
(1-\eta)\eta^k\Omega(k)
=
g_u(k)\Omega(k).
\]
Finally, $(1-\eta)\eta^k>0$ and $\pi_{k+1}\ge0$, so
$\Omega(k)\ge0$ for every $k\in\Nzero$.
\end{proof}

\begin{theorem}
	\label{lem:kernel}
	Let $u>0$, and let $H$ be a real polynomial of degree at most
	four satisfying $H(k)\ge0$ for every $k\in\Nzero$. With
	\[
	g_u(k)=\frac1{1+u}\left(\frac{u}{1+u}\right)^k,
	\]
	we have
	\begin{equation}\label{eq:kernel}
		\Psi_u(H):=
		\sum_{k=0}^{\infty}(k-7)(k-14)H(k)g_u(k)\ge0.
	\end{equation}
	If $H\ne0$, equality is possible only when $u=1$. At that
	parameter, equality holds precisely for
	\begin{equation}\label{eq:kerneleq}
		H(x)=
		\alpha x(x-1)(x-2)(x-3)
		+\beta x(x-1)(x-3)(x-4),
		\qquad \alpha,\beta\ge0.
	\end{equation}
\end{theorem}
The series defining $\Psi_u(H)$ converges absolutely. Indeed,
with $\eta=u/(1+u)\in(0,1)$, its summands satisfy
\[
\bigl|(k-7)(k-14)H(k)g_u(k)\bigr|
\le C(1+k^6)\eta^k
\]
for a constant $C$ depending on $H$ and $u$, and the series on
the right converges. The inequality and its equality
characterization are proved in Appendix~\ref{app:kernel},

\begin{proof}[Proof of part (i) of of Theorem~\ref{thm:main} ]
Consider
\[
\phi_*(z)=
\left(\frac{1+z/10}{1-z/10}\right)^5.
\]
Its shifted generating function is
\[
F_*(w)=\phi_*(10(w-1))
=\frac{w^5}{(2-w)^5}
=\sum_{j=0}^{\infty}
\binom{j+4}{4}\frac{w^{j+5}}{2^{j+5}}.
\]
The series has nonnegative coefficients and converges for
$|w|<2$. Hence every derivative of $F_*$ is nonnegative on
$[0,1]$, and the change of variables $w=1+z/10$ proves
absolute monotonicity of $\phi_*$ on $[-10,0]$.
Direct differentiation gives
$F_*(1)=1, F_*'(1)=10, F_*''(1)=100$.
Thus, $\phi_*(0)=\phi_*'(0)=\phi_*''(0)=1$, so $\phi_*\in\mathcal B_{5,2}$.
Furthermore,
\[
F_*(w)=\frac1{32}w^5+O(w^6).
\]
Thus by Lemma~\ref{lem:gap}, we have $\rad(\phi_*)=10$. Let $\phi\in\mathcal B_{5,2}$ and suppose that $r>0$ is
feasible, meaning that $\phi$ is absolutely monotonic on
$[-r,0]$. Write
\[
F(w)=\phi(r(w-1))
=\sum_{j=0}^{\infty}\pi_jw^j.
\]
Absolute monotonicity gives $\pi_j\ge0$, and
\[
\sum_{j=0}^{\infty}\pi_j=F(1)=\phi(0)=1.
\]
Thus the coefficients define a random variable $X_F$ by
$\mathbb P (X_F=j)=\pi_j$. The second-order conditions give
\[
\mathbb E[X_F]=F'(1)=r,\qquad
\mathbb E[X_F(X_F-1)]=F''(1)=r^2,
\]
and hence
\[
\mathbb E[X_F^2]=r^2+r.
\]
By Lemma~\ref{lem:s5tail},
\[
\pi_{k+1}=g_u(k)\Omega(k),
\qquad
u=\frac rb>0,
\]
where $\Omega$ has degree at most four and is nonnegative on
$\Nzero$. Define
\[
q(x)=(x-8)(x-15)=x^2-23x+120.
\]
Since $q(0)=120$ and $q(k+1)=(k-7)(k-14)$,
Lemma~\ref{lem:kernel} yields
\[
\begin{aligned}
\mathbb E[q(X_F)]
&=120\pi_0+
\sum_{k=0}^{\infty}(k-7)(k-14)\pi_{k+1}\\
&=120\pi_0+\Psi_u(\Omega)\ge0.
\end{aligned}
\]
On the other hand, the moment identities imply
\begin{equation}\label{eq:forbidden}
\begin{aligned}
\mathbb E[q(X_F)]
&=\mathbb E[X_F^2]-23\mathbb E[X_F]+120\\
&=r^2+r-23r+120\\
&=(r-10)(r-12).
\end{aligned}
\end{equation}
The last expression is negative for $10<r<12$. Therefore no feasible radius belongs to $(10,12)$.
This also excludes every radius above twelve. Indeed, if
$\rad(\phi)>10$, the definition of the supremum supplies a
feasible radius $R>10$. Choose
\[
10<r<\min\{R,12\}.
\]
Since $[-r,0]\subseteq[-R,0]$, the radius $r$ is also feasible,
contradicting \eqref{eq:forbidden}. Thus
\[
\rad(\phi)\le10
\qquad\text{for every }\phi\in\mathcal B_{5,2}.
\]
Together with the attaining function $\rad(\phi_*)=10$, prove that $B_{5,2} = 10$.
Now we show that $\phi_*$ is the only function with radius of absolute monotonicity equal $10$.   
Suppose that $\phi$ is feasible at $r=10$.
Equation~\eqref{eq:forbidden} gives
\[
0=120\pi_0+\Psi_u(\Omega).
\]
Both terms are nonnegative, so
\[
\pi_0=0,\qquad \Psi_u(\Omega)=0.
\]
Moreover, $\Omega\ne0$, since otherwise every probability at
a positive integer would vanish, contradicting $\mathbb E[X_F]=10$.
Lemma~\ref{lem:kernel} therefore implies
\[
u=1,\qquad
\Omega=\alpha h_0+\beta h_1,
\qquad \alpha,\beta\ge0,
\]
where
\[
h_0(x)=x(x-1)(x-2)(x-3),
\qquad
h_1(x)=x(x-1)(x-3)(x-4).
\]
Let $Y$ have the geometric distribution
\[
\Pr(Y=k)=g_1(k)=2^{-k-1},
\qquad k\in\Nzero.
\]
Its generating function is $(2-w)^{-1}$, so its factorial moments satisfy
\[
\mathbb E[Y^{\underline j}]=j!,
\qquad
\mathbb E[(Y+1)Y^{\underline j}]=2(j+1)!,
\qquad j\ge0.
\]
The second identity follows from
\[
(Y+1)Y^{\underline j}=
Y^{\underline{j+1}}+(j+1)Y^{\underline j}.
\]	
Using
\[
h_0(x)=x^{\underline4},
\qquad
h_1(x)=x^{\underline4}-2x^{\underline3}
+2x^{\underline2},
\]
we obtain
\[
	\begin{array}{c|cc}
		&\mathbb E[h_i(Y)]&\mathbb E[(Y+1)h_i(Y)]\\ \hline
		h_0&24&240\\
		h_1&16&168
	\end{array}
	\]
For any function $f$, we have 
\begin{equation}
\label{exp12}
\mathbb{E}[f(X_F)]
=
\pi_0 f(0)
+
\mathbb{E}\!\left[f(Y+1)\Omega(Y)\right],
\end{equation}
whenever the expectations are absolutely convergent. Indeed, this is
simply the identity
\[
\sum_{j=0}^{\infty} f(j)\pi_j
=
\pi_0 f(0)
+
\sum_{k=0}^{\infty} f(k+1)g_1(k)\Omega(k).
\]
Taking into account that $\pi_0=0$ and using $f(x)=1$ and $f(x)=x$ in (\ref{exp12}) give	
\[
\mathbb E[\Omega(Y)]=1,
\qquad
\mathbb E[(Y+1)\Omega(Y)]=10.
\]
Consequently,
\[
24\alpha+16\beta=1,
\qquad
240\alpha+168\beta=10.
\]
Solving for $\alpha$ and $\beta$, we obtain 	
\[
\beta=0,\qquad \alpha=\frac1{24},
\qquad \Omega(x)=\frac{x^{\underline4}}{24}.
\]
The generating function is therefore
\[
F(w)
=w\sum_{k=0}^{\infty}g_1(k)\Omega(k)w^k
=\frac{w}{24}\,
w^4\frac{d^4}{dw^4}\left(\frac1{2-w}\right)
=\frac{w^5}{(2-w)^5}.
\]
Since $u=r/b=1$ and $r=10$, we also have $b=10$. Returning
to the original variable $z$ yields
\[
\phi(z)=F(1+z/10)
=\left(\frac{1+z/10}{1-z/10}\right)^5
=\phi_*(z).
\]
Finally, suppose that $\rad(\phi)=10$ is defined only as a
supremum. Every radius smaller than ten is feasible by
nestedness. The denominator $(1-z/b)^5$, with $b>0$, has no
zero on the nonpositive real axis, so every derivative of
$\phi$ is continuous there. Taking feasible radii tending to
10 gives, for each $j\ge0$,
\[
\phi^{(j)}(-10)
=\lim_{r\uparrow10}\phi^{(j)}(-r)\ge0.
\]
Thus $\phi$ is feasible at radius ten, and the preceding
equality argument applies. This proves uniqueness of the
maximizer.
\end{proof}

\begin{remark}
The new scalar counterexamples do not contradict the stage-level result of Kocsis and N\'emeth \cite[Theorem~2.6]{KN}, which bounds second-order diagonally implicit SSP Runge--Kutta coefficients by $2s$. In particular, \eqref{eq:sevenphi} is not A-stable, and none of the present constructions is asserted to realize its scalar radius as a nonlinear SSP coefficient.
\end{remark}

\section{Conclusions}
In this work, for second order, we establish the conjectured optimum
$B_{5,2}=10$, extending the cases $s=1,2,3,4$ proved in
\cite{LK}. However, the family \eqref{eq:secondfamily}
disproves the conjectured value $B_{s,2}=2s$ for every
$s\ge6$. At third order, the conjectured value holds for
$s=2,3,4$ by \cite{LK}, whereas the family
\eqref{eq:thirdF} disproves it for every $s\ge5$.

The exact value of $B_{5,3}$, the remaining optimal values
at higher multiplicities, and the sharp asymptotic growth
of $B_{s,p}$ remain to be determined. The product construction
provides an asymptotic lower bound, but does not establish
the growth rate of the optimum.

Our results concern scalar radii of absolute monotonicity.
They do not contradict the stage-level result of Kocsis and
N\'emeth \cite[Theorem~2.6]{KN}, which bounds the SSP
coefficient of every second-order, $s$-stage diagonally
implicit Runge--Kutta method by $2s$. We do not claim that
any of our scalar counterexamples can be realized by a
Runge--Kutta method whose nonlinear SSP coefficient equals
its scalar radius.

\appendix
\section{Proof of Theorem \ref{lem:kernel}}\label{app:kernel}
In this appendix, we prove Lemma~\ref{lem:kernel}. The proof is
divided into two parts. For $0<u\le3/2$, we establish
nonnegativity by expressing $\Psi_u(H)$ as an explicit
nonnegative linear combination of values of $H$ at nonnegative
integers. For $u>3/2$, we derive an integral representation of
$\Psi_u(H)$ and establish strict positivity for nonzero $H$
by proving that the associated moment matrices are positive
definite.
\subsection{An explicit positive representation of $\Psi_u$ for $0<u\le3/2$}

\begin{proposition}
\label{lem:five-node}
For every real polynomial $H$ of degree at most four, the functional $\Psi_u(H)$ in (\ref{eq:kernel})
is given by
\begin{equation}\label{eq:five-node}
\Psi_u(H)
=A_0(u)H(0)+A_1(u)H(1)+A_2(u)H(2)
+A_3(u)H(3)+A_{15}(u)H(15),
\end{equation}
where
\begin{align*}
A_0(u)
&=
2(2u^2-7u+7)(2u^4-2u^3+5u^2+7),\\
A_1(u)
&=
-\frac{6u}{7}
\left(
30u^5-130u^4+254u^3-273u^2+182u-91
\right),\\
A_2(u)
&=
\frac{60}{13}u^2(u-1)^2(6u^2-13u+13),\\
A_3(u)
&=
-2u^3(5u^3-20u^2+33u-22),\\
A_{15}(u)
&=
\frac{2}{91}u^4(u-1)^2.
\end{align*}
For $0<u\le3/2$, these weights satisfy
\[
A_0(u)>0,\qquad A_1(u)>0,\qquad A_3(u)>0,
\qquad A_2(u)\ge0,\qquad A_{15}(u)\ge0.
\]
The last two weights vanish precisely when $u=1$.
Consequently, if $H$ is nonnegative on $\mathbb Z_{\ge0}$, then
\[
\Psi_u(H)\ge0 \quad \text{for} \quad 
0<u\le3/2.
\]
For $H\ne0$, this inequality is strict whenever $u\ne1$.
At $u=1$, the representation reduces to
\begin{equation}\label{eq:five-node-u-one}
\Psi_1(H)=48H(0)+24H(1)+8H(3).
\end{equation}
At that parameter, equality holds precisely for
\begin{equation}\label{eq:kerneleq}
	H(x)=
	\alpha x(x-1)(x-2)(x-3)
	+\beta x(x-1)(x-3)(x-4),
	\qquad \alpha,\beta\ge0.
\end{equation}
\end{proposition}
\begin{proof}
The generating function of the weights $g_{u}(k)$ is
\[
\sum_{k=0}^{\infty}g_u(k)w^k
=
\frac{1}{1+u-uw},
\qquad
|w|<\frac{1+u}{u}.
\]
Since $(1+u)/u>1$, this power series may be differentiated
term by term in a neighborhood of $w=1$. Its $j$th derivative
is
\[
\sum_{k=0}^{\infty}k^{\underline j}g_u(k)w^{k-j}
=
\frac{j!u^j}{(1+u-uw)^{j+1}}.
\]
Evaluating at $w=1$ gives
\begin{equation}\label{eq:factorial-moments}
\sum_{k=0}^{\infty}k^{\underline j}g_u(k)=j!u^j,
\qquad j\ge0.
\end{equation}
All sums involving polynomial factors converge absolutely, because the weights decay geometrically.
Using
\[
(x-7)(x-14)x^{\underline j}
=x^{\underline{j+2}} +(2j-20)x^{\underline{j+1}}
+(j-7)(j-14)x^{\underline j},
\]
and applying \eqref{eq:factorial-moments}, we obtain
\begin{equation}\label{eq:kernel-factorial-moments}
m_j(u):= \Psi_u(x^{\underline j})=
j!u^j\Bigl[(j+1)(j+2)u^2+(2j-20)(j+1)u
+(j-7)(j-14)
\Bigr].
\end{equation}
In particular,
\begin{equation}\label{mexpression}
\begin{split}
m_0(u)&=2(u^2-10u+49),\quad  m_1(u)=6u(u^2-6u+13),\\
m_2(u)&=24u^2(u^2-4u+5), \quad  m_3(u)=24u^3(5u^2-14u+11),\\
m_4(u)&=720u^4(u-1)^2.
\end{split}
\end{equation}
Let $\mathcal N=\{0,1,2,3,15\}$, and define the Lagrange polynomials
\[
\ell_a(x)
=
\prod_{\substack{b\in\mathcal N\\b\ne a}}
\frac{x-b}{a-b},
\qquad a\in\mathcal N.
\]
They satisfy
\[
\ell_a(b)=
\begin{cases}
1,&b=a,\\
0,&b\ne a,
\end{cases}
\qquad a,b\in\mathcal N.
\]
Since $\deg H\le4$, interpolation gives the exact identity
\[
H(x)=\sum_{a\in\mathcal N}H(a)\ell_a(x).
\]
Applying the linear functional $\Psi_u$ therefore yields
$\Psi_u(H) = \sum_{a\in\mathcal N}\Psi_u(\ell_a)H(a)$.
Thus the required weights in (\ref{eq:five-node}) are
$A_a(u)=\Psi_u(\ell_a)$. To compute them explicitly, we expand the Lagrange polynomials in the falling-factorial basis:
\begin{align*}
\ell_0(x) &=
1-x^{\underline1}
+\frac12x^{\underline2}
-\frac16x^{\underline3}
+\frac1{90}x^{\underline4},\quad 
\ell_1(x)
=
x^{\underline1}
-x^{\underline2}
+\frac12x^{\underline3}
-\frac1{28}x^{\underline4},\\
\ell_2(x)
&=
\frac12x^{\underline2}
-\frac12x^{\underline3}
+\frac1{26}x^{\underline4},\quad
\ell_3(x)
=
\frac16x^{\underline3}
-\frac1{72}x^{\underline4},\\
\ell_{15}(x)
&=
\frac1{32760}x^{\underline4}.
\end{align*}
Consequently,
\begin{align*}
A_0(u)
&=
m_0(u)-m_1(u)+\frac12m_2(u)
-\frac16m_3(u)+\frac1{90}m_4(u),\\
A_1(u)
&=
m_1(u)-m_2(u)+\frac12m_3(u)-\frac1{28}m_4(u),\\
A_2(u)
&=
\frac12m_2(u)-\frac12m_3(u)+\frac1{26}m_4(u),\\
A_3(u)
&=
\frac16m_3(u)-\frac1{72}m_4(u),\\
A_{15}(u)
&=
\frac1{32760}m_4(u).
\end{align*}
Substituting the expressions from (\ref{mexpression}) and collecting
terms gives the five weights stated in the lemma.
This proves \eqref{eq:five-node} for every $u>0$.
Now we show the positivity of the weights for $0<u\le3/2$. First,
\[
2u^2-7u+7 =2\left(u-\frac74\right)^2+\frac78>0,
\]
and
\[
2u^4-2u^3+5u^2+7 =
2u^2\left(u-\frac12\right)^2+\frac92u^2+7>0.
\]
Hence $A_0(u)>0$.	
The quadratic $6u^2-13u+13$ has positive leading coefficient
and negative discriminant. Therefore $A_2(u)\ge0$, with equality precisely at $u=1$.
The same conclusion for $A_{15}(u)$ follows immediately from its formula.
For the other two weights, define
\[
p(u)=30u^5-130u^4+254u^3-273u^2+182u-91
\]
and
\[
q(u)=5u^3-20u^2+33u-22.
\]
We show that both polynomials are strictly increasing.
For $q$, we have
\[
q'(u)
=
15\left(u-\frac43\right)^2+\frac{19}{3}>0.
\]
For $p$, writing $u=1+t$ gives
\[
\frac32p'(1+t)
=t^2(15t+4)^2+137t^2+27t+42.
\]
The quadratic on the right is strictly positive, since its
leading coefficient is positive and its descriminant is negative.
Thus $p'(u)>0$ for every real $u$. Since
$p(3/2)=-\frac{85}{16}<0$ and $q(3/2)=-\frac58<0$,
monotonicity implies
$p(u)<0,\qquad q(u)<0$ for $0<u\le3/2$.
The formulas for $A_1$ and $A_3$ now show that 
$A_1(u)>0$ and $A_3(u)>0$ for $0<u\le3/2$.
Thus, we conclude that  if $H$ is nonnegative on $\mathbb Z_{\ge0}$ then
\[
\Psi_u(H)\ge0
\qquad(0<u\le3/2).
\]
When $u\ne1$, all five weights are strictly positive.
Equality would force
\[
H(0)=H(1)=H(2)=H(3)=H(15)=0.
\]
A polynomial of degree at most four with five distinct
zeros is identically zero. Hence the inequality is strict
for nonzero $H$ when $u\ne1$. Finally, at $u=1$ the weights are
\[
A_0(1)=48,\qquad A_1(1)=24,\qquad A_3(1)=8,
\qquad A_2(1)=A_{15}(1)=0.
\]
Substitution into \eqref{eq:five-node} proves
\eqref{eq:five-node-u-one}. 
The preceding arguments show that equality for nonzero $H$ requires
$u=1$. By \eqref{eq:five-node-u-one}, it is then equivalent to
$H(0)=H(1)=H(3)=0$. Consequently
\[
H(x)=x(x-1)(x-3)(ax+b)
\]
for real $a,b$. Nonnegativity at $x=2$ and $x=4$ gives
\[
2a+b\le0,\qquad4a+b\ge0.
\]
Set
\[
\alpha=\frac{4a+b}{2}\ge0,
\qquad
\beta=-\frac{2a+b}{2}\ge0.
\]
Then
\[
H(x)=\alpha x(x-1)(x-2)(x-3)
+\beta x(x-1)(x-3)(x-4).
\]
Conversely, these two generators are nonnegative on $\Nzero$ and vanish
at $0,1,3$, so every such combination gives equality at $u=1$.
\end{proof}

\subsection{An integral representation of $\Psi_u$}
\begin{proposition}
Let $H$ be a polynomial of degree at most $4$ non-negative in $\mathbb Z_{\ge0}$ written in 
the the falling-factorial basis as
\[
H(x)=\sum_{j=0}^{4}c_jx^{\underline j}.
\]
Define
\[
\widehat H(y)=\sum_{j=0}^{4}c_ju^jy^j.
\]
Then $\widehat H \not=0$ and $\widehat H(y) \geq 0$ for $y \geq 0$. Moreover, we have 
\begin{equation}\label{eq:kernel-large-u-integral}
\Psi_u(H)
=\int_0^\infty
e^{-y}
\bigl(v^2y^2-24vy+120\bigr)\widehat H(y)\,dy \quad \text{where} \quad  v = 1+u.
\end{equation}
\end{proposition}
\begin{proof}
For every integer $j\ge0$ and every $t\ge0$,
\[
e^{-t}
\sum_{k=0}^{\infty}
k^{\underline j}\frac{t^k}{k!}
=
e^{-t}t^j
\sum_{k=j}^{\infty}\frac{t^{k-j}}{(k-j)!}
=
t^j.
\]
Consequently,
\begin{equation}\label{eq:kernel-positive-transform}
\widehat H(y)=
e^{-uy}
\sum_{k=0}^{\infty}H(k)\frac{(uy)^k}{k!},
\qquad y\ge0.
\end{equation}
The series converges absolutely, since $H$ is a polynomial.
Every term on the right-hand side is nonnegative, so
\[
\widehat H(y)\ge0 \quad \text{for} \quad 
y\ge0.
\]
Moreover, $\widehat H\ne0$ since $u>0$ and thus the coefficients
$c_ju^j$ can all vanish only if all $c_j$ vanish.
To show (\ref{eq:kernel-large-u-integral}), we notice that both sides depend linearly on $H$, so it suffices to verify
the identity for $H(x)=x^{\underline j}$, $0\le j\le4$.
From (\ref{eq:kernel-factorial-moments}), we have 
\[
\Psi_u(x^{\underline j})
=
j!u^j
\bigl(
(j+1)(j+2)v^2-24(j+1)v+120
\bigr), \quad v= u+1. 
\]
On the other hand, the transform of $x^{\underline j}$ is
$u^jy^j$, and
\[
\begin{aligned}
	\int_0^\infty
	e^{-y}
	\bigl(v^2y^2-24vy+120\bigr)u^jy^j\,dy & =
	u^j\bigl(v^2(j+2)!-24v(j+1)!+120j!\bigr)\\
	&=
	j!u^j
	\bigl(
	(j+1)(j+2)v^2-24(j+1)v+120
	\bigr).
\end{aligned}
\]
This proves \eqref{eq:kernel-large-u-integral}.	
\end{proof} 
We shall need the following Lemma.
\begin{lemma}
Every real polynomial $Q$ of degree at most four that is
nonnegative on $[0,\infty)$ admits a representation
\begin{equation}\label{eq:kernel-half-line-squares}
Q(y)=
\sum_{\nu=1}^{2}a_\nu(y)^2
+
y\sum_{\nu=1}^{2}b_\nu(y)^2,
\qquad
\deg a_\nu\le2,\quad
\deg b_\nu\le1.
\end{equation}
Some of the polynomials in this representation may be zero.
\end{lemma}
\begin{proof}
Observe that $Q(t^2)$ is a globally nonnegative real polynomial 
of degree at most eight. Its real zeros have even multiplicity, and its nonreal zeros
occur in conjugate pairs. Factoring over $\mathbb C$ and
assigning one zero from each conjugate pair to a complex
polynomial therefore gives
\[
Q(t^2)=U(t)^2+V(t)^2
\]
for real polynomials $U,V$ of degree at most four.
Separate their even and odd parts:
\[
U(t)=a_1(t^2)+t\,b_1(t^2),
\qquad
V(t)=a_2(t^2)+t\,b_2(t^2),
\]
where $\deg a_\nu\le2$ and $\deg b_\nu\le1$.
Since $Q(t^2)$ is even, averaging the expressions at $t$ and
$-t$ eliminates the cross terms and gives
\[
Q(t^2)
=
a_1(t^2)^2+a_2(t^2)^2
+t^2\bigl(b_1(t^2)^2+b_2(t^2)^2\bigr).
\]
Replacing $t^2$ by $y$ proves
\eqref{eq:kernel-half-line-squares}.
\end{proof}

\begin{proposition}
\label{lem:kernel-large-u}
Let $H$ be a polynomial of degree at most $4$ non-negative in $\mathbb Z_{\ge0}$. 
For $u\ge3/2$,
\[
\Psi_u(H)\ge0,
\]
with equality if and only if $H$ is the zero polynomial.
\end{proposition}
\begin{proof}
Define the linear functional
\[
\mathcal L_v[f]
=
\int_0^\infty
e^{-y}\bigl(v^2y^2-24vy+120\bigr)f(y)\,dy.
\]
Its moments are given by
\[
\mu_j(v)
=
\mathcal L_v[y^j]
=
j!\bigl((j+1)(j+2)v^2-24(j+1)v+120\bigr).
\]
Introduce the symmetric matrices
\[
M(v)
=
\begin{pmatrix}
\mu_0(v)&\mu_1(v)&\mu_2(v)\\
\mu_1(v)&\mu_2(v)&\mu_3(v)\\
\mu_2(v)&\mu_3(v)&\mu_4(v)
\end{pmatrix}
\]
and
\[
N(v)
=
\begin{pmatrix}
\mu_1(v)&\mu_2(v)\\
\mu_2(v)&\mu_3(v)
\end{pmatrix}.
\]
For
\[
a(y)=a_0+a_1y+a_2y^2,
\qquad
b(y)=b_0+b_1y,
\]
we have
\[
\mathcal L_v[a^2]
=
\begin{pmatrix}a_0&a_1&a_2\end{pmatrix}
M(v)
\begin{pmatrix}a_0\\a_1\\a_2\end{pmatrix}
\]
and
\[
\mathcal L_v[yb^2]
=
\begin{pmatrix}b_0&b_1\end{pmatrix}
N(v)
\begin{pmatrix}b_0\\b_1\end{pmatrix}.
\]
We shall first prove positive definiteness of $M(v)$ and $N(v)$ for 
$v\ge5/2$. By Sylvester criterion, a real symmetric matrix is positive
definite if and only if all its leading principal minors
are strictly positive. For $M(v)$, direct determinant expansion gives the leading
principal minors
\[
2\bigl((v-6)^2+24\bigr),
\qquad
12A(v),
\qquad
576B(v),
\]
where
\begin{align*}
A(v)
&=
v^4-24v^3+256v^2-960v+1200,\\
B(v)
&=
v^6-36v^5+588v^4-4176v^3
+14040v^2-21600v+12000.
\end{align*}
The first minor is strictly positive for every real $v$.	
For $v\ge3$, write $v=3+t$, where $t\ge0$.
The identities
\[
A(3+t)
=
t^2(t-6)^2+58t^2+36t+57
\]
and
\[
B(3+t)=
t^4(t-9)^2+
102\left(t^2-\frac32\right)^2
+9t^2+180t^3+270t+\frac{375}{2}
\]
show that $A(v)>0$ and $B(v)>0$.
For $5/2\le v\le3$, write $v=3-t$, where
$0\le t\le1/2$. Then
\[
A(3-t)
=t^4+12t^3+94t^2-36t+57
\ge57-36\left(\frac12\right)
=39>0.
\]
Similarly,
\[
B(3-t)
=
t^6+18t^5+183t^4-180t^3-297t^2-270t+417.
\]
Discarding its nonnegative terms and using $t\le1/2$ gives
\[
B(3-t)\ge
417-180\left(\frac12\right)^3
-297\left(\frac12\right)^2
-270\left(\frac12\right)=
\frac{741}{4}>0.
\]
Thus every leading principal minor of $M(v)$ is strictly
positive for $v\ge5/2$. Sylvester's criterion proves that
$M(v)$ is positive definite throughout this range.
The leading principal minors of $N(v)$ are
\[
6\bigl((v-4)^2+4\bigr) \quad \text{and} \quad 	144(v-2)C(v),
\]
where
\[
C(v)=v^3-14v^2+70v-100.
\]
The first minor is strictly positive. Moreover,
\[
C'(v)
=
3\left(v-\frac{14}{3}\right)^2+\frac{14}{3}>0,
\]
so $C$ is strictly increasing. Since
\[
C\left(\frac52\right)=\frac{25}{8}>0,
\]
we have $C(v)>0$ for every $v\ge5/2$.
Also $v-2>0$ in this range. Hence the second minor is
strictly positive, and Sylvester's criterion proves that
$N(v)$ is positive definite.

Now, apply \eqref{eq:kernel-half-line-squares} to $\widehat H$.
Using \eqref{eq:kernel-large-u-integral}, we obtain
\[
\Psi_u(H)
=
\sum_{\nu=1}^{2}\mathcal L_v[a_\nu^2]
+
\sum_{\nu=1}^{2}\mathcal L_v[yb_\nu^2], \quad v = u+1 \geq 5/2.
\]
The positive definiteness of $M(v)$ and $N(v)$ shows that
every summand is nonnegative, and that a summand is strictly
positive whenever its corresponding polynomial is nonzero.
Since $\widehat H\ne0$, at least one of the polynomials
$a_\nu,b_\nu$ is nonzero. Consequently,
\[
\Psi_u(H)>0.
\]
This proves the proposition.
\end{proof}
Combining  Propositions \ref{lem:five-node} and \ref{lem:kernel-large-u} 
proves Theorem \ref{lem:kernel}.
	

\bigskip 

{\bf{Declarations}} The author used AI tools to assist with manuscript drafting and selected computational tasks. The author has verified all mathematical arguments and takes full responsibility for the content of this work.

\end{document}